\documentclass[leqno,draft]{article}

\def\diam{\mathop{\rm diam}}

\begin{document}

\title{Notes on $p$-adic numbers}

\author{Stephen William Semmes	\\
	Rice University		\\
	Houston, Texas}

\date{}

\maketitle

	As usual the real numbers are denoted ${\bf R}$.  For each $x
\in {\bf R}$, the \emph{absolute value} of $x$ is denoted $|x|$ and
defined to be $x$ when $x \ge 0$ and to be $-x$ when $x \le 0$.  Thus
$|0| = 0$ and $|x| > 0$ when $x \ne 0$.  One can check that $|x + y|
\le |x| + |y|$ and $|x \, y| = |x| \, |y|$ for every $x, y \in {\bf
R}$.

	Let $N(x)$ be a nonnegative real-valued function defined on
the rational numbers ${\bf Q}$ such that $N(0) = 0$, $N(x) > 0$ when
$x \ne 0$,
\begin{equation}
\label{N(x y) = N(x) N(y)}
	N(x \, y) = N(x) \, N(y)
\end{equation}
for all $x, y \in {\bf Q}$, and
\begin{equation}
\label{N(x + y) le C (N(x) + N(y))}
	N(x + y) \le C \, (N(x) + N(y))
\end{equation}
for some $C \ge 1$ and all $x, y \in {\bf Q}$.  For the usual
\emph{triangle inequality} one aks that this condition holds with $C =
1$, i.e.,
\begin{equation}
\label{N(x + y) le N(x) + N(y)}
	N(x + y) \le N(x) + N(y)
\end{equation}
for all $x, y \in {\bf Q}$.  The \emph{ultrametric} version of the
triangle inequality is stronger still and asks that
\begin{equation}
\label{N(x + y) le max (N(x), N(y))}
	N(x + y) \le \max (N(x), N(y))
\end{equation}
for all $x, y \in {\bf Q}$.  If $N$ satisfies (\ref{N(x + y) le C
(N(x) + N(y))}), $l$ is a positive integer, and $x_1, \ldots, x_{2^l}
\in {\bf Q}$, then
\begin{equation}
\label{N(sum_{j=1}^{2^l} x_j) le C^l sum_{j=1}^{2^l} N(x_j)}
	N\Big(\sum_{j=1}^{2^l} x_j\Big)
		\le C^l \, \sum_{j=1}^{2^l} N(x_j),
\end{equation}
as one can check using induction on $l$.

	The usual absolute value function $|x|$ satisfies these
conditions with the ordinary triangle inequality (\ref{N(x + y) le max
(N(x), N(y))}).  If $N(x) = 0$ when $x = 0$ and $N(x) = 1$ when $x \ne
0$, then $N(x)$ satisfies these conditions with the ultrametric
version of the triangle inequality.  For each prime number $p$, the
\emph{$p$-adic absolute value} of a rational number $x$ is denoted
$|x|_p$ and defined by $|x|_p = 0$ when $x = 0$, and $|x|_p = p^{-j}$
when $x$ is equal to $p^j$ times a ratio of nonzero integers, neither
of which is divisible by $p$.  One can check that $|x|_p$ satisfies
these conditions with the ultrametric version of the triangle
inequality.  Roughly speaking, this says that $x + y$ has at least $j$
factors of $p$ when $x$ and $y$ have at least $j$ factors of $p$.

	Let $a$ be a positive real number.  Clearly
\begin{equation}
	\max (r, t) \le (r^a + t^a)^{1/a}
\end{equation}
for all $r, t \ge 0$, and hence
\begin{equation}
	r + t \le \max(r^{1-a}, t^{1-a}) \, (r^a + t^a)
		\le (r^a + t^a)^{1/a}
\end{equation}
when $0 < a \le 1$, which implies that
\begin{equation}
	(r + t)^a \le r^a + t^a.
\end{equation}
When $a \ge 1$,
\begin{equation}
	(r + t)^a \le 2^{a - 1} \, (r^a + t^a),
\end{equation}
because $x^a$ is a convex function on the nonnegative real numbers.

	For all $a > 0$, $N(x)^a$ is a nonnegative real-valued
function on ${\bf Q}$ which vanishes at $0$, is positive at all
nonzero $x \in {\bf Q}$, and sends products to products.  If $N(x)$
satisfies (\ref{N(x + y) le C (N(x) + N(y))}), then
\begin{equation}
	N(x + y)^a \le C^a \, (N(x)^a + N(y)^a)
\end{equation}
when $0 < a \le 1$ and $N(x + y)^a \le 2^{a-1} \, C^a \, (N(x)^a +
N(y)^a)$ when $a \ge 1$.  In particular, if $N(x)$ satisfies the
ordinary triangle inequality (\ref{N(x + y) le N(x) + N(y)}) and $0 <
a \le 1$, then $N(x)^a$ also satisfies the ordinary triangle
inequality.  If $N(x)$ satisfies the ultrametric version (\ref{N(x +
y) le max (N(x), N(y))}) of the triangle inequality, then $N(x)^a$
satisfies the ultrametric version of the triangle inequality for all
$a > 0$.

	Observe that
\begin{equation}
	N(1) = 1,
\end{equation}
because $N(1)^2 = N(1)$ and $N(1) > 0$, and similarly
\begin{equation}
	N(-1) = 1.
\end{equation}
Suppose that $N(x) > 1$ for some integer $x$.  In this event $N$ is
unbounded on the integers ${\bf Z}$, because
\begin{equation}
\label{N(x^n) = N(x)^n to infty}
	N(x^n) = N(x)^n \to \infty
\end{equation}
as $n \to \infty$, and we say that $N$ is \emph{Archimedian}.
Otherwise,
\begin{equation}
	N(x) \le 1 \quad\hbox{for all } x \in {\bf Z},
\end{equation}
and $N$ is \emph{non-Archimedian}.

	Suppose that $N$ is Archimedian, and let $A$ be the set of
positive real numbers $a$ such that $N(x)$ is bounded by a constant
times $|x|^a$ for every $x \in {\bf Z}$.  It follows from
(\ref{N(sum_{j=1}^{2^l} x_j) le C^l sum_{j=1}^{2^l} N(x_j)}) that $A
\ne \emptyset$, and (\ref{N(x^n) = N(x)^n to infty}) implies that $A$
has a positive lower bound.  For every $a \in A$ and $x \in {\bf Z}$,
$N(x) \le |x|^a$, because $N(x) > |x|^a$ implies that $N(x^j)$ is not
bounded by a constant times $|x|^{j \, a}$.  Let $\alpha$ be the
infimum of $A$.  By the previous remarks, $\alpha > 0$ and
\begin{equation}
	N(x) \le |x|^\alpha
\end{equation}
for every $x \in {\bf Z}$.

	Let $n$ be an integer, $n > 1$.  Every nonnegative integer can
be expressed as a finite sum of the form $\sum_{j=0}^l r_j \, n^j$,
where each $r_j$ is an integer and $0 \le r_j < n$.  If $N(n) \le 1$,
then $N(x)$ is bounded by a constant times a power of the logarithm of
$1 + |x|$ for $x \in {\bf Z}$, contradicting the fact that $A$ has a
positive lower bound.  Similarly, if $N(n) < n^\alpha$, which means
that $N(n) = n^b$ for some $b < \alpha$, then $N(x)$ is bounded by a
constant times a power of the logarithm of $1 + |x|$ times $|x|^b$ for
$x \in {\bf Z}$, contradicting the fact that $\alpha$ is the infimum
of $A$.  We conclude that
\begin{equation}
	N(x) = |x|^\alpha
\end{equation}
for every $x \in {\bf Z}$, and therefore for every $x \in {\bf Q}$.

	Now suppose that $N$ is non-Archimedian, and let $x, y \in
{\bf Q}$ be given.  The \emph{binomial theorem} asserts that
\begin{equation}
	(x + y)^n = \sum_{j=0}^n {n \choose j} \, x^j \, y^{n-j}
\end{equation}
for each positive integer $n$, where the \emph{binomial coefficients}
\begin{equation}
	{n \choose j} = \frac{n!}{j! \, (n - j)!}
\end{equation}
are integers, and $k!$ is ``$k$ factorial'', the product of the
positive integers from $1$ to $k$, which is interpreted as being equal
to $1$ when $k = 0$.  If $l$ is a positive integer and $n + 1 \le
2^l$, then (\ref{N(sum_{j=1}^{2^l} x_j) le C^l sum_{j=1}^{2^l}
N(x_j)}) implies that
\begin{eqnarray}
	N(x + y)^n & = & N\big((x+y)^n\big) 			\\
		& \le & C^l \sum_{j=0}^n N(x)^j \, N(y)^{n-j}	\nonumber \\
		& \le & (n + 1) \, C^l \max (N(x)^n, N(y)^n).	\nonumber
\end{eqnarray}
Using this one can check that $N$ satisfies the ultrametric version of
the triangle inequality (\ref{N(x + y) le max (N(x), N(y))}).

	If $N(x) = 1$ for each nonzero integer $x$, then $N(x) = 1$
for every $x \in {\bf Q}$, $x \ne 0$.  Otherwise, there is an integer
$p > 1$ such that $N(p) < 1$, and we may choose $p$ to be as small as
possible.  If $p = p_1 \, p_2$, where $p_1, p_2 \in {\bf Z}$ and $p_1,
p_2 > 1$, then $N(p_1) \, N(p_2) < 1$, and hence $N(p_1) < 1$ or
$N(p_2) < 1$, which contradicts the minimality of $p$.  Therefore $p$
is prime.  Of course $N(x) \le 1$ for every $x \in {\bf Z}$, and the
minimality of $p$ implies that $N(x) = 1$ when $x \in {\bf Z}$ and $1
\le x < p$.  If $x = p \, y$, $y \in {\bf Z}$, then
\begin{equation}
	N(x) = N(p) \, N(y) \le N(p) < 1.
\end{equation}
Using the ultrametric version of the triangle inequality for $N$, one
can check that $N(x) = 1$ when $x \in {\bf Z}$ and $x$ is not an
integer multiple of $p$.  Moreover,
\begin{equation}
	N(x) = |x|_p^a
\end{equation}
for every $x \in {\bf Q}$, where $a > 0$ is determined by $N(p) = p^{-a}$.

	Fix a prime number $p$.  The \emph{$p$-adic metric} on ${\bf
Q}$ is defined by
\begin{equation}
	d_p(x, y) = |x - y|_p
\end{equation}
for $x, y \in {\bf Q}$.  A sequence $\{x_j\}_{j=1}^\infty$ of rational
numbers \emph{converges} to $x \in {\bf Q}$ in the $p$-adic metric if
for every $\epsilon > 0$ there is an $L \ge 1$ such that
\begin{equation}
	|x_j - x|_p < \epsilon
\end{equation}
for every $j \ge L$.  If $\{x_j\}_{j=1}^\infty$,
$\{y_j\}_{j=1}^\infty$ are sequences of rational numbers which
converge in the $p$-adic metric to $x, y \in {\bf Q}$, then the
sequences of sums $x_j + y_j$ and products $x_j \, y_j$ converge to
the sum $x + y$ and product $x \, y$ of the limits of the initial
sequences, by standard arguments.

	A sequence $\{x_j\}_{j=1}^\infty$ of rational numbers
is a \emph{Cauchy sequence} with respect to the $p$-adic metric
if for each $\epsilon > 0$ there is an $L \ge 1$ such that
\begin{equation}
	|x_j - x_l|_p < \epsilon
\end{equation}
for every $j, l \ge L$.  Every convergent sequence in ${\bf Q}$ is a
Cauchy sequence.  If $\{x_j\}_{j=1}^\infty$ is a Cauchy sequence in
${\bf Q}$ with respect to the $p$-adic metric, then the sequence of
differences $x_j - x_{j+1}$ converges to $0$ in the $p$-adic metric.
Of course the analogous statement also works for the standard metric
$|x - y|$.  For the $p$-adic metric, the converse holds because of the
ultrametric version of the triangle inequality.

	The \emph{$p$-adic numbers} ${\bf Q}_p$ are obtained by
completing the rational numbers with respect to the $p$-adic metric,
just as the real numbers are obtained by completing the rational
numbers with respect to the standard metric.  To be more precise,
${\bf Q} \subseteq {\bf Q}_p$, and ${\bf Q}_p$ is equipped with
operations of addition and multiplication which agree with the usual
arithmetic operations on ${\bf Q}$ and which satisfy the standard
field axioms.  There is an extension of the $p$-adic absolute value
$|x|_p$ to $x \in {\bf Q}_p$ which satisfies the same conditions as
before, i.e., $|x|_p$ is a nonnegative real number for all $x \in {\bf
Q}_p$ which is equal to $0$ exactly when $x = 0$, and
\begin{equation}
	|x \, y|_p = |x|_p \, |y|_p
\end{equation}
and
\begin{equation}
	|x + y|_p \le \max (|x|_p, |y|_p)
\end{equation}
for every $x, y \in {\bf Q}_p$.  We can extend the $p$-adic metric to
${\bf Q}_p$ using the same formula
\begin{equation}
	d_p(x, y) = |x - y|_p
\end{equation}
for $x, y \in {\bf Q}_p$.  The rational numbers are dense in ${\bf
Q}_p$, in the sense that for each $x \in {\bf Q}_p$ and $\epsilon > 0$
there is a $y \in {\bf Q}$ such that
\begin{equation}
	|x - y|_p < \epsilon.
\end{equation}
A sequence $\{x_j\}_{j=1}^\infty$ of $p$-adic numbers converges to $x
\in {\bf Q}_p$ if for each $\epsilon > 0$ there is an $L \ge 1$ such
that
\begin{equation}
	|x_j - x|_p < \epsilon
\end{equation}
for every $j \ge L$.  One can show that sums and products of
convergent sequences in ${\bf Q}_p$ converge to the sums and products
of the limits of the individual sequences.  A sequence
$\{x_j\}_{j=1}^\infty$ in ${\bf Q}_p$ is a Cauchy sequence if
for every $\epsilon > 0$ there is an $L \ge 1$ such that
\begin{equation}
	|x_j - x_l|_p < \epsilon
\end{equation}
for every $j, l \ge L$.  Convergent sequences are Cauchy sequences,
and a sequence $\{x_j\}_{j=1}^\infty$ in ${\bf Q}_p$ is a Cauchy
sequence if and only if
\begin{equation}
	\lim_{j \to \infty} (x_j - x_{j+1}) = 0 \quad\hbox{in } {\bf Q}_p,
\end{equation}
because of the ultrametric version of the triangle inequality.  The
$p$-adic numbers are \emph{complete} in the sense that every Cauchy
sequence in ${\bf Q}_p$ converges to an element of ${\bf Q}_p$.

	Suppose that $x \in {\bf Q}_p$, $x \ne 0$.  Because ${\bf Q}$
is dense in ${\bf Q}_p$, there is a $y \in {\bf Q}$ such that
\begin{equation}
	|x - y|_p < |x|_p.
\end{equation}
Using the ultrametric version of the triangle inequality, one can
check that
\begin{equation}
	|x|_p = |y|_p.
\end{equation}
It follows that $|x|_p$ is an integer power of $p$.

	An infinite series $\sum_{j=0}^\infty a_j$ with terms in ${\bf
Q}_p$ \emph{converges} if the corresponding sequence of partial sums
$\sum_{j=0}^n a_j$ converges in ${\bf Q}_p$.  It is easy to see that
the partial sums for $\sum_{j=0}^\infty a_j$ form a Cauchy sequence in
${\bf Q}_p$ if and only if the $a_j$'s converge to $0$ in ${\bf Q}_p$,
in which event the series converges.  If $\sum_{j=0}^\infty a_j$,
$\sum_{j=0}^\infty b_j$ are convergent series in ${\bf Q}_p$ and
$\alpha, \beta \in {\bf Q}_p$, then $\sum_{j=0}^\infty (\alpha \, a_j
+ \beta \, b_j)$ converges, and
\begin{equation}
	\sum_{j=0}^\infty (\alpha \, a_j + \beta \, b_j)
	  = \alpha \sum_{j=0}^\infty a_j + \beta \, \sum_{j=0}^\infty b_j.
\end{equation}
This follows from the analogous statements for convergent sequences.

	Suppose that $x \in {\bf Q}_p$.  For each nonnegative integer
$n$,
\begin{equation}
	(1 - x) \sum_{j=0}^n x^j = 1 - x^{n+1}.
\end{equation}
Here $x^j$ is interpreted as being $1$ when $j = 0$ for every $x \in
{\bf Q}_p$.  If $x \ne 1$, then we can rewrite this identity as
\begin{equation}
	\sum_{j=0}^n x^j = \frac{1 - x^{n+1}}{1 - x}.
\end{equation}
If $|x|_p < 1$, then $\sum_{j=0}^\infty x^j$ converges, and
\begin{equation}
\label{sum_{j=0}^infty x^j = frac{1}{1 - x}}
	\sum_{j=0}^\infty x^j = \frac{1}{1 - x}.
\end{equation}

	Let us write $p \, {\bf Z}$ for the set of integer multiples
of $p$, which is an ideal in the ring of integers ${\bf Z}$.  There is
a natural mapping from ${\bf Z}$ onto ${\bf Z} / p \, {\bf Z}$, the
integers modulo $p$.  Addition and multiplication are well-defined
modulo $p$, and this mapping is a homomorphism, which is to say that
it preserves these operations.  The fact that $p$ is prime implies
that ${\bf Z} / p \, {\bf Z}$ is a field, which means that nonzero
elements have inverses.  Equivalently, for each integer $x$ which is
not a multiple of $p$, there is an integer $y$ such that $x \, y - 1
\in p \, {\bf Z}$.

	For every $w \in {\bf Z}$, $\sum_{j=0}^\infty p^j \, w^j$
converges in ${\bf Q}_p$ to $1 / (1 - p \, w)$, as in
(\ref{sum_{j=0}^infty x^j = frac{1}{1 - x}}).  In particular, $1 / (1
- p \, w)$ is the limit of a sequence of integers in the $p$-adic
metric.  Hence if $a, b \in {\bf Z}$ and $b - 1 \in p \, {\bf Z}$,
then $a / b$ is the limit of a sequence of integers in the $p$-adic
metric.  If $b \in {\bf Z}$ and $b \not\in p \, {\bf Z}$, then there
is a $c \in {\bf Z}$ such that $b \, c - 1 \in p \, {\bf Z}$, and
hence $a / b + (a \, c) / (b \, c)$ is a limit of integers in the
$p$-adic metric.

	Of course every integer has $p$-adic absolute value less than
or equal to $1$, and the limit of any sequence of integers in the
$p$-adic metric has $p$-adic absolute value less than or equal to $1$.
The discussion in the previous paragraph shows that every $x \in {\bf
Q}$ with $|x|_p \le 1$ is the limit of a sequence of integers in the
$p$-adic metric.  In other words,
\begin{equation}
	\{x \in {\bf Q} : |x|_p \le 1\}
\end{equation}
is the same as the closure of ${\bf Z}$ in ${\bf Q}$ with respect to
the $p$-adic metric.

	Put
\begin{equation}
	{\bf Z}_p = \{x \in {\bf Q}_p : |x|_p \le 1\}.
\end{equation}
Every element of ${\bf Q}_p$ is the limit of a sequence of rational
numbers in the $p$-adic metric, because ${\bf Q}$ is dense in ${\bf
Q}_p$.  If $x \in {\bf Z}_p$, then $x$ is the limit of a sequence of
rational numbers $\{x_j\}_{j=1}^\infty$ in the $p$-adic metric, and
$|x_j|_p \le 1$ for sufficiently large $j$ by the ultrametric version
of the triangle inequality.  Because rational numbers with $p$-adic
absolute value less than or equal to $1$ can be approximated by
integers in the $p$-adic metric, $x$ is the limit of a sequence of
integers in the $p$-adic metric.  In short, ${\bf Z}_p$ is equal to
the closure of ${\bf Z}$ in ${\bf Q}_p$.

	The elements of ${\bf Z}_p$ are said to be \emph{$p$-adic
integers}.  Let us write $p \, {\bf Z}_p$ for the set of $p$-adic
integers of the form $p \, x$, $x \in {\bf Z}_p$, which is the same as
\begin{equation}
	p \, {\bf Z}_p = \{y \in {\bf Q}_p : |y|_p \le 1/p\}.
\end{equation}
Of course ${\bf Z} \subseteq {\bf Z}_p$ and $p \, {\bf Z} \subseteq p
\, {\bf Z}_p$.  If $x \in {\bf Z}_p$, then $x$ can be expressed as $y
+ w$, where $y \in {\bf Z}$ and $w \in p \, {\bf Z}_p$.

	Note that sums and products of $p$-adic integers are $p$-adic
integers, which is to say that ${\bf Z}_p$ is a ring with respect to
addition and multiplication, a subring of the field ${\bf Q}_p$.
Furthermore, $p \, {\bf Z}_p$ is an ideal in ${\bf Z}_p$, which means
that sums of elements of $p \, {\bf Z}_p$ lie in $p \, {\bf Z}_p$, and
the product of an element of $p \, {\bf Z}_p$ and an element of ${\bf
Z}_p$ lies in $p \, {\bf Z}_p$.  Consequently, there is a quotient
ring ${\bf Z}_p / p \, {\bf Z}_p$ and a natural mapping from ${\bf
Z}_p$ onto the quotient ${\bf Z}_p / p \, {\bf Z}_p$ which preserves
addition and multiplication.

	The inclusions of $p \, {\bf Z}$, ${\bf Z}$ in $p \, {\bf
Z}_p$, ${\bf Z}_p$, respectively, lead to a natural mapping
\begin{equation}
\label{{bf Z} / p {bf Z} to {bf Z}_p / p {bf Z}_p}
	{\bf Z} / p \, {\bf Z} \to {\bf Z}_p / p \, {\bf Z}_p.
\end{equation}
If one maps an element of ${\bf Z}$ into ${\bf Z} / p \, {\bf Z}$, and
then into ${\bf Z}_p / p \, {\bf Z}_p$, then that is the same as
mapping the element of ${\bf Z}$ into ${\bf Z}_p$, and then into the
quotient ${\bf Z}_p / p \, {\bf Z}_p$.  One can check that the mapping
(\ref{{bf Z} / p {bf Z} to {bf Z}_p / p {bf Z}_p}) is a ring
isomorphism.  To describe the inverse more explicitly, if $x \in {\bf
Z}_p$, $x = y + w$, $y \in {\bf Z}$, $w \in p \, {\bf Z}_p$, then the
image of $y$ in the quotient ${\bf Z} / p \, {\bf Z}$ corresponds
exactly to the image of $x$ in the quotient ${\bf Z}_p / p \, {\bf
Z}_p$.

	More generally, for each positive integer $j$, $p^j \, {\bf
Z}$ denotes the set of integer multiples of $p^j$, which is an ideal
in ${\bf Z}$, and $p^j \, {\bf Z}_p$ denotes the ideal in ${\bf Z}_p$
consisting of the products $p^j \, x$, $x \in {\bf Z}_p$, which is the
same as
\begin{equation}
	p^j \, {\bf Z}_p = \{y \in {\bf Q}_p : |y|_p \le p^{-j}\}.
\end{equation}
The inclusion $p^j \, {\bf Z} \subseteq p^j \, {\bf Z}_p$ leads to a
homomorphism
\begin{equation}
	{\bf Z} / p^j \, {\bf Z} \to {\bf Z}_p / p^j \, {\bf Z}_p
\end{equation}
between the corresponding quotients, which is an isomorphism.  This
uses the fact that each $x \in {\bf Z}_p$ can be expressed as $y + w$,
with $y \in {\bf Z}$ and $w \in p^j \, {\bf Z}_p$, because ${\bf Z}$
is dense in ${\bf Z}_p$.

	We say that $\mathcal{C} \subseteq {\bf Q}_p$ is a \emph{cell}
if it is of the form
\begin{equation}
	\mathcal{C} = \{y \in {\bf Q}_p : |y - x|_p \le p^l\}
\end{equation}
for some $x \in {\bf Q}_p$ and $l \in {\bf Z}$.  The diameter of
$\mathcal{C}$, denoted $\diam \mathcal{C}$, is equal to $p^{-l}$ in
this case, which is to say that the maximum of the $p$-adic distances
between elements of $\mathcal{C}$ is equal to $p^{-l}$.  Note that
cells are both open and closed as subsets of ${\bf Q}_p$, and hence
${\bf Q}_p$ is totally disconnected.  More precisely, the distances
between points in a cell and points in the complement are greater than
or equal to the diameter of the cell.  Also, if $\mathcal{C}_1$,
$\mathcal{C}_2$ are cells in ${\bf Q}_p$, then $\mathcal{C}_1
\subseteq \mathcal{C}_2$, or $\mathcal{C}_2 \subseteq \mathcal{C}_1$,
or $\mathcal{C}_1 \cap \mathcal{C}_2 = \emptyset$.

	If $\mathcal{C}$ is a cell in ${\bf Q}_p$ and $n$ is a
positive integer, then $\mathcal{C}$ contains $p^n$ disjoint cells
with diameter equal to $p^{-n} \, \mathcal{C}$, and $\mathcal{C}$ is
equal to the union of these smaller cells.  For instance, ${\bf Z}_p$
is equal to the disjoint union of translates of $p^n \, {\bf Z}_p$,
one for each element of ${\bf Z} / p^n \, {\bf Z}$.  Every cell in
${\bf Q}_p$ can be obtained from ${\bf Z}_p$ by a translation and
dilation, and so the decompositions for ${\bf Z}_p$ yield the
analogous decompositions for arbitrary cells in ${\bf Q}_p$.  One can
show that cells in ${\bf Q}_p$ are compact, in practically the same
way as for closed and bounded intervals in the real line.
Consequently, closed and bounded subsets of ${\bf Q}_p$ are compact.

	Suppose that $\mathcal{C}$ is a cell in ${\bf Q}_p$ and that
$f(x)$ is a continuous real-valued function on $\mathcal{C}$.  By
standard results in analysis, $f$ is uniformly continuous.  For each
positive integer $n$, let $\mathcal{C}_{1, n}, \ldots,
\mathcal{C}_{p^n, n}$ be the $p^n$ disjoint cells of diameter $p^{-n}
\, \diam \mathcal{C}$ contained in $\mathcal{C}$.  This leads to the
Riemann sums
\begin{equation}
	\sum_{j=1}^{p^n} f(x_j) \, p^{-n} \, \diam \mathcal{C},
\end{equation}
where $x_j \in \mathcal{C}_{j, n}$ for each $j$.  Because of uniform
continuity, these Riemann sums converge as $n \to \infty$ to the
\emph{integral of $f$ over $\mathcal{C}$}, and the limit does not
depend on the choices of the $x_j$'s.

	Of course a ${\bf Q}_p$-valued function on $\mathcal{C}$ is
uniformly continuous too, since $\mathcal{C}$ is compact.  Riemann
sums for ${\bf Q}_p$-valued functions can be defined in the same way
as in the previous paragraph, since $p^{-n} \, \diam \mathcal{C} \in
{\bf Q}$.  However,
\begin{equation}
	|p^{-n}|_p = p^n \to \infty \hbox{ as } n \to \infty,
\end{equation}
and uniform continuity is not sufficient to imply convergence of the
Riemann sums as $n \to \infty$.  As in Section 12.4 in \cite{C}, one
can look at ${\bf Q}_p$-valued measures on $\mathcal{C}$ which are
bounded on the cells contained in $\mathcal{C}$.  Because of the
ultrametric version of the triangle inequality, boundedness of the
measures of small cells is adequate to get convergence of Riemann sums
of continuous ${\bf Q}_p$-valued functions on $\mathcal{C}$.

	One can also consider continuous ${\bf Q}_\ell$-valued
functions on $\mathcal{C}$, where $\ell$ is prime and $\ell \ne p$.
As in the previous situations, such a function is uniformly continuous
because of the compactness of $\mathcal{C}$.  Integer powers of $p$
have absolute value equal to $1$ in ${\bf Q}_\ell$, which implies that
the absolute value of a Riemann sum of a ${\bf Q}_\ell$-valued
function on $\mathcal{C}$ is bounded by the maximum of the absolute
value of the function on $\mathcal{C}$, because of the ultrametric
version of the triangle inequality.  Uniform continuity of the
function implies convergence of the Riemann sums, because differences
of Riemann sums can be bounded by maximal local oscillations of the
function.

	Let $f(x)$ be a polynomial on ${\bf Q}_p$, which means that
\begin{equation}
	f(x) = a_n \, x^n + \cdots + a_0
\end{equation}
for some nonnegative integer $n$ and $a_0, \ldots, a_n \in {\bf Q}_p$.
For the discussion that follows it is convenient to ask that the
coefficients $a_j$ of $f(x)$ be $p$-adic integers.

	For each $x, y \in {\bf Z}_p$,
\begin{equation}
\label{|f(x) - f(y)|_p le |x - y|_p}
	|f(x) - f(y)|_p \le |x - y|_p.
\end{equation}
To see this one can first consider the case where $y = 0$, and then
use the observation that the translate of a polynomial with $p$-adic
integer coefficients by an element of ${\bf Z}_p$ also has $p$-adic
integer coefficients.  Alternatively, one can use the identity
\begin{equation}
	x^n - y^n = (x - y) \, \sum_{j=0}^{n-1} x^j \, y^{n-1-j}
\end{equation}
for $n \ge 1$.

	The derivative of $f(x)$ is the polynomial
\begin{equation}
	f'(x) = n \, a_n \, x^{n-1} + \cdots + a_1.
\end{equation}
Thus $f'(x)$ has $p$-adic integer coefficients, and therefore
\begin{equation}
\label{|f'(x) - f'(y)|_p le |x - y|_p}
	|f'(x) - f'(y)|_p \le |x - y|_p
\end{equation}
for every $x, y \in {\bf Z}_p$.

	Observe that
\begin{equation}
\label{|f(x) - f(y) - f'(y) (x - y)|_p le |x - y|_p^2}
	|f(x) - f(y) - f'(y) \, (x - y)|_p \le |x - y|_p^2
\end{equation}
for every $x, y \in {\bf Z}_p$.  One can check this first when $y =
0$, and then use translations to get the general case.  One can also
show this when $f(x) = x^n$ and then sum over $n$.

	Suppose that $z \in {\bf Q}_p$, $f(z) \in p \, {\bf Z}_p$, and
that $|f'(z)|_p = 1$.  \emph{Hensel's lemma} asserts that there is a
$w \in {\bf Z}_p$ such that $w - z \in p \, {\bf Z}_p$ and $f(w) = 0$.

	To prove this we use Newton's method.  Put $x_0 = z$, and for
each $j \ge 1$ choose $x_j$ according to the rule
\begin{equation}
	f(x_{j-1}) + f'(x_{j-1}) \, (x_j - x_{j-1}) = 0.
\end{equation}
Equivalently,
\begin{equation}
	x_j = x_{j-1} + \frac{f(x_{j-1})}{f'(x_{j-1})}.
\end{equation}
More precisely,
\begin{equation}
	x_{j-1} \in {\bf Z}_p, \quad f(x_{j-1}) \in p \, {\bf Z}_p,
		\quad |f'(x_{j-1})|_p = 1
\end{equation}
by induction, which implies that
\begin{equation}
	x_j - x_{j-1} \in p \, {\bf Z}_p.
\end{equation}
This ensures that the analogous conditions hold for $x_j$, by
(\ref{|f(x) - f(y)|_p le |x - y|_p}) and (\ref{|f'(x) - f'(y)|_p le |x
- y|_p}).

	Using (\ref{|f(x) - f(y) - f'(y) (x - y)|_p le |x - y|_p^2}),
we get
\begin{equation}
	|f(x_j)|_p \le |x_j - x_{j-1}|_p^2.
\end{equation}
Since $|x_{j+1} - x_j|_p \le |f(x_j)|_p$,
\begin{equation}
	|x_{j+1} - x_j|_p \le |x_j - x_{j-1}|_p^2.
\end{equation}
It follows that $\lim_{j \to \infty} (x_{j+1} - x_j) = 0$ in ${\bf
Q}_p$, and hence that $\{x_j\}_{j=1}^\infty$ converges to an element
$w$ of ${\bf Q}_p$.  Because $x_j - x_{j-1} \in p \, {\bf Z}_p$ for
every $j$, $w - z \in p \, {\bf Z}_p$ and $w \in {\bf Z}_p$.
Moreover, $\lim_{j \to \infty} f(x_j) = 0$ in ${\bf Q}_p$, and
therefore $f(w) = 0$, as desired.

	Suppose now that
\begin{equation}
	z \in {\bf Z}_p, \quad |f(z)|_p < |f'(z)|_p^2.
\end{equation}
If $f(z) \in p \, {\bf Z}_p$ and $|f'(z)|_p = 1$, then $|f(z)|_p <
|f'(z)|_p^2$, and thus these conditions are more general than the
previous ones.  A refined version of Hensel's lemma asserts that there
is a $w \in {\bf Z}_p$ such that $f(w) = 0$.

	Again we put $x_0 = z$ and choose $x_j$ according to the same
rule as in the previous situation.  This uses the induction hypotheses
\begin{equation}
	x_{j-1} \in {\bf Z}_p, 
		\quad |f(x_{j-1})|_p < |f'(x_{j-1})|_p^2 = |f'(z)|_p^2.
\end{equation}
Under these conditions,
\begin{equation}
	|x_j - x_{j-1}|_p \le \frac{|f(x_{j-1})|_p}{|f'(x_{j-1})|_p}
				< |f'(x_{j-1})|_p,
\end{equation}
which implies that $|f'(x_j) - f'(x_{j-1})|_p < |f'(x_{j-1})|_p$ and
$|f'(x_j)|_p = |f'(x_{j-1})|_p$.  Furthermore,
\begin{equation}
	|f(x_j)|_p \le |x_j - x_{j-1}|_p^2 < |f(x_{j-1})|_p,
\end{equation}
and hence the induction hypotheses continue to be satisfied for $x_j$.
One can check that
\begin{equation}
	\lim_{j \to \infty} f(x_j) = 0
\end{equation}
in ${\bf Q}_p$, and that $\{x_j\}_{j=1}^\infty$ converges to $w \in
{\bf Z}_p$ such that $f(w) = 0$.

	If $x, y \in {\bf Q}_p$, $n$ is a positive integer, and $x =
y^n$, then $|x|_p = |y|_p^n$, which implies that $x = 0$ or $|x|_p$ is
of the form $p^{l \, n}$ for some $l \in {\bf Z}$.  Of course $x =
y^n$ with $y = 0$ when $x = 0$.  If $|x|_p = p^{l \, n}$ for some
integer $l$, and $x_1 = p^{-l} \, x$, then $|x_1|_p = 1$ and $x = y^n$
for some $y \in {\bf Q}_p$ if and only if $x_1 = (y_1)^n$ for some
$y_1 \in {\bf Q}_p$ with $|y_1|_p = 1$.

	The mapping
\begin{equation}
	h(y) = y^n
\end{equation}
satisfies $h({\bf Z}_p) \subseteq {\bf Z}_p$ and $|h(z) - h(w)|_p \le
|z - w|_p$.  Therefore $h$ induces a mapping on ${\bf Z}_p / p \, {\bf
Z}_p \cong {\bf Z} / p \, {\bf Z}$, which also takes $n$th powers.

	Fix a positive integer $q$ which is prime and $a \in {\bf
Q}_p$ with $|a|_p = 1$.  Consider the polynomial
\begin{equation}
	f(x) = x^q - a.
\end{equation}
Thus the coefficients of $f$ are $p$-adic integers, and the zeros of
$f$ are the $q$th roots of $a$.  The derivative of $f$ is
\begin{equation}
	f'(x) = q \, x^{q-1}.
\end{equation}
In particular, for every $x \in {\bf Q}_p$ with $|x|_p = 1$,
\begin{equation}
	|f'(x)|_p = 1
\end{equation}
when $q \ne p$ and
\begin{equation}
	|f'(x)|_p = \frac{1}{p}
\end{equation}
when $q = p$.

	Let $({\bf Z} / p \, {\bf Z})^*$ be the group of nonzero
elements of ${\bf Z} / p \, {\bf Z}$ under multiplication, a finite
abelian group with $p - 1$ elements.  It is a well-known theorem that
every finite abelian group is isomorphic to a Cartesian product of
finitely many cyclic groups.  Another well-known theorem asserts that
${\bf Z} / p \, {\bf Z}$ is cyclic.

	If $p - 1$ is not an integer multiple of $q$, then every
element of ${\bf Z} / p \, {\bf Z}$ is a $q$th power.  If $p - 1$ is
not an integer multiple of $q$ and $q \ne p$, then we can apply
Hensel's lemma to get that $f(x) = x^q - a$ has a root.

	Suppose that $p - 1$ is an integer multiple of $q$, which
implies that $q \ne p$.  A necessary condition for $f(x) = x^q - a$ to
have a root in ${\bf Z}_p$ is that the image of $a$ in ${\bf Z}_p / p
\, {\bf Z}_p \cong {\bf Z} / p \, {\bf Z}$ be a $q$th power there.
Hensel's lemma implies that this condition is sufficient too.

	When $q = p = 2$, one can check that $x^2 - 1 \in 8 \, {\bf
Z}_2$ for every $x \in {\bf Z}_2$ such that $|x|_2 = 1$.  Conversely,
the refined version of Hensel's lemma implies that every $y \in {\bf
Z}_2$ such that $y - 1 \in 8 \, {\bf Z}_2$ is equal to $x^2$ for some
$x \in {\bf Q}_2$ such that $|x|_2 = 1$.

	Let us consider real and complex-valued functions again
briefly.  Suppose that $f(x)$ is a continuous function from ${\bf
Z}_p$ into ${\bf R}$ or ${\bf C}$ such that
\begin{equation}
	f(x + y) = f(x) + f(y)
\end{equation}
for every $x, y \in {\bf Z}_p$.  The image of $f$ is a compact
subgroup of ${\bf R}$ or ${\bf C}$, as appropriate, under addition,
and it follows that $f(x) = 0$ for every $x \in {\bf Z}_p$.

	Now suppose that $\phi(x)$ is a continuous function from ${\bf
Z}_p$ into nonzero complex numbers such that
\begin{equation}
	\phi(x + y) = \phi(x) \, \phi(y)
\end{equation}
for every $x, y \in {\bf Z}_p$.  Since $\log |\phi(x)|$ is a
homomorphism from ${\bf Z}_p$ into ${\bf R}$ with respect to addition,
the remarks of the previous paragraph yield
\begin{equation}
	|\phi(x)| = 1
\end{equation}
for every $x \in {\bf Z}_p$.  If $U$ is a small neighborhood of $1$ in
${\bf C}$, then there is a positive integer $j$ such that
\begin{equation}
	\phi(p^j \, {\bf Z}_p) \subseteq U,
\end{equation}
because of the continuity of $\phi$.  Since $p^j \, {\bf Z}_p$ is a
subgroup of ${\bf Z}_p$ under addition, $\phi(p^j \, {\bf Z}_p)$ is a
subgroup of the unit circle in ${\bf C}$ under multiplication, and
it follows that
\begin{equation}
	\phi(x) = 1
\end{equation}
for every $x \in p^j \, {\bf Z}_p$.  Thus $\phi$ corresponds to a
homomorphism from ${\bf Z}_p / p^j \, {\bf Z}_p \cong {\bf Z} / p^j \,
{\bf Z}$ as a group under addition into the unit circle in ${\bf C}$
as a group under multiplication.

	Next let us consider some connections with matrices, following
\cite{C}.  Of course $1 \times 1$ matrices are scalars, and we start
with them.

	Let $p$ be a prime and let $x$ be an element of ${\bf Z}_p$
such that $x \ne 1$ and $x - 1 \, p \, {\bf Z}_p$.  We can express $x$
as $1 + p^j \, a$, where $j$ is a positive integer and $|a|_p = 1$.
Let $q$ be a prime and consider $x^q$.

	Using the binomial theorem we can express $x^q$ as $1 + q \,
p^j \, a + b \, p^{2 \, j}$, where $b \in {\bf Z}_p$.  If $q \ne p$ or
$j \ge 2$, then $x^q - 1$ reduces to $q \, p^j \, a \ne 0$ in ${\bf
Z}_p / p^{j+1} \, {\bf Z}_p$, and therefore $x^q \ne 1$.

	Suppose that $q = p$, $j = 1$, and $p \ne 2$.  Using the
binomial theorem we can express $x^p$ as $1 + p^2 \, a + p^3 \, b$,
where $b \in {\bf Z}_p$.  Again we get $x^p \ne 1$.

	When $q = p = 2$ and $j = 1$, $x^2 = 1$ is possible, since $x$
may be $-1$.  As oberseved previously, $|x|_2 = 1$ implies that $x^2 -
1 \in 8 \, {\bf Z}_2$.

	To summarize, when $p \ne 2$, $x^q \ne 1$, and it follows that
$x^t \ne 1$ for all positive integers $t$, by repeating the argument.
When $p = 2$, it may be that $x^2 = 1$.  Otherwise, $x^t \ne 1$ for
all positive integers $t$.

	Fix a positive integer $n$, and let $\mathcal{M}_n({\bf Q}_p)$
be the space of $n \times n$ matrices with entries in ${\bf Q}_p$.  We
can add and multiply matrices in $\mathcal{M}_n({\bf Q}_p)$ in the
usual way, or multiply matrices and elements of ${\bf Q}_p$.

	Let $\mathcal{M}_n({\bf Z}_p)$ be the space of $n \times n$
matrices with entries in ${\bf Z}_p$.  Sums and products of matrices
in $\mathcal{M}_n({\bf Z}_p)$ also lie in $\mathcal{M}_n({\bf Z}_p)$,
as do products of matrices in $\mathcal{M}_n({\bf Z}_p)$ and elements
of ${\bf Z}_p$.  Of course $I \in \mathcal{M}_n({\bf Z}_p)$, where the
entries of $I$ on the diagonal are equal to $1$ and the entries off of
the diagonal are equal to $0$.

	For each $A \mathcal{M}_n({\bf Q}_p)$ the \emph{determinant}
$\det A$ is defined in the usual way and is an element of ${\bf Q}_p$.
If $A \in \mathcal{M}_n({\bf Z}_p)$, then $\det A \in {\bf Z}_p$.

	We say that $A \in \mathcal{M}_n({\bf Q}_p)$ is
\emph{invertible} if there is a $B \in \mathcal{M}_n({\bf Q}_p)$ such
that $A \, B = B \, A = I$, in which event the inverse $B$ of $A$ is
denoted $A^{-1}$.  Standard results in linear algebra imply that $A
\in \mathcal{M}_n({\bf Q}_p)$ is invertible if and only if $\det A \ne
0$.  If $A \in \mathcal{M}_n({\bf Z}_p)$, then $A$ is invertible and
$A^{-1} \in \mathcal{M}_n({\bf Z}_p)$ if and only if $|\det A|_p = 1$.

	Suppose that $A \in \mathcal{M}_n({\bf Z}_p)$, $A \ne I$, and
$A - I$ has entries in $p \, {\bf Z}_p$.  Hence $A = I + p^j \, B$ for
some positive integer $j$ and $B \in \mathcal{M}_n({\bf Z}_p)$, where
at least one entry of $B$ has $p$-adic absolute value equal to $1$.

	Let $q$ be prime and consider $A^q$.  The binomial theorem
implies that
\begin{equation}
	A^q - I - q \, p^j \, B
\end{equation}
has entries in $p^{2 \, j} \, {\bf Z}_p$, and it follows that $A^q \ne
I$ when $q \ne p$ or $j \ge 2$.  If $q = p \ne 2$ and $j = 1$, then
\begin{equation}
	A^p - I - p^2 \, B
\end{equation}
has entries in $p^3 \, {\bf Z}_p$, and $A^p \ne I$.

	If $q = p = 2$ and $j = 1$, then
\begin{equation}
	A^2 = I + 4 \, B + 4 \, B^2.
\end{equation}
It may be that $A^2 = I$, and in any event $A^2 - I$ has entries in $4
\, {\bf Z}_2$, and the previous discussion applies to $A^2$ if $A^2
\ne I$.

	Let $GL_n({\bf Q}_p)$ be the group of $n \times n$ invertible
matrices with entries in ${\bf Q}_p$, and let $GL_n({\bf Z}_p)$ be the
subgroup of $GL({\bf Q}_p)$ consisting of matrices with entries in
${\bf Z}_p$ whose inverse have the same property.  There is a natural
homomorphism from $GL_n({\bf Z}_p)$ onto $GL_n({\bf Z} / p \, {\bf
Z})$, the group of $n \times n$ invertible matrices with entries in
the field ${\bf Z} / p \, {\bf Z}$, induced by the ring homomorphism
from ${\bf Z}_p$ into ${\bf Z}_p / p \, {\bf Z}_p \cong {\bf Z} / p \,
{\bf Z}$.  Of course the determinant defines a homomorphism from
$GL_n({\bf Q}_p)$ into the group of nonzero elements of ${\bf Q}_p$
under multiplication, and from $GL_n({\bf Z}_p)$ into the group of
$p$-adic numbers with $p$-adic absolute value equal to $1$ under
multiplication.

	If $p \ne 2$, $A \in GL_n({\bf Z}_p)$, $A - I$ has entries in
$p \, {\bf Z}_p$, and $A \ne I$, then $A^q \ne I$ for every prime $q$,
and therefore $A^l \ne I$ for every positive integer $l$.  If $G$ is a
subgroup of $GL_n({\bf Z}_p)$ with only finitely many elements, then
for each $A \in G$ there is a positive integer $l$ such that $A^l =
I$, and it follows that the natural homomorphism from $GL_n({\bf
Z}_p)$ into $GL_n({\bf Z} / p \, {\bf Z})$ is one-to-one on $G$,
because $I$ is the only element of $G$ which can be sent to the
identity in $GL_n({\bf Z} / p \, {\bf Z})$.

	If $p = 2$, $A \in GL_n({\bf Z}_2)$, $A - I$ has entries in $2
\, {\bf Z}_2$, and $A \ne I$, then $A^2 = I$ or $A^l \ne I$ for every
positive integer $l$.  If $H$ is a subgroup of $GL_n({\bf Z}_2)$ such
that $H$ has only finitely many elements and every entry of $A - I$ is
an element of $p \, {\bf Z}_2$ for every $A \in H$, then $A^2 = I$ for
every $A \in H$, and hence $H$ is abelian, since $(A \, B)^2 = A \, B
\, A \, B = I$ implies that $A \, B = B \, A$ when $A^2 = B^2 = I$.
The earlier remarks also imply that $A = I$ when $A \in H$ and $A - I$
has entries in $4 \, {\bf Z}_2$.

	If $G$ is a subgroup of $GL_n({\bf Z}_2)$ with only finitely
many elements, then the subgroup $H$ of $G$ consisting of the matrices
$A$ such that $A - I$ has entries in $2 \, {\bf Z}_2$ satisfies the
conditions described in the previous paragraph.  Of course $H$ is the
same as the subgroup of $A \in G$ which go to the identity under the
natural homomorphism from $GL_n({\bf Z}_2)$ to $GL_n({\bf Z} / 2 \,
{\bf Z})$.

	Let $GL_n({\bf Q})$ be the group of $n \times n$ invertible
matrices with entries in ${\bf Q}$.  We can think of $GL_n({\bf Q})$
as a subgroup of $GL_n({\bf Q}_p)$ for every prime $p$.  If $G$ is a
subgroup of $GL_n({\bf Q})$ with only finitely many elements, then $G$
is actually a subgroup of $GL_n({\bf Z}_p)$ for all but finitely many
$p$.

	If $H \subseteq GL_n({\bf Q}_2)$ is the collection of diagonal
matrices with diagonal entries $\pm 1$, then $H$ is a subgroup of
$GL_n({\bf Z}_2)$ and every $A \in H$ satisfies $A^2 = I$ and $A - I$
has entries in $2 \, {\bf Z}_2$.  The conjugate of $H$ by any element
of $GL_2({\bf Z}_2)$ has the same features.

	If $A$ is a linear transformation on a vector space $V$ over a
field $k$ with characteristic $\ne 2$ such that $A^2 = I$, then $V$ is
spanned by the eigenspaces of $V$ corresponding to the eigenvalues
$\pm 1$.  Specifically,
\begin{equation}
	P_1 = \frac{I - A}{2}, \quad P_2 = \frac{I + A}{2}
\end{equation}
are the projections of $V$ onto these eigenspaces, and any linear
transformation $B$ on $V$ which commutes with $A$ maps these
eigenspaces onto themselves.  These facts provide additional
information about subgroups $H$ of $GL_n({\bf Z}_2)$ such that $A^2 =
I$ for every $A \in H$.

\end{document}